%% file: rzamalgam.tex
\documentclass[a4paper,10pt]{amsart}
 \usepackage{amsfonts,mathrsfs}
 \usepackage{amsthm}
 \usepackage{amssymb}
 \usepackage{enumerate}
 \usepackage{tikz-cd}
 \usepackage{cleveref}
 \usepackage{pgf,tikz}
\usetikzlibrary{arrows}
\theoremstyle{definition}
\newtheorem{Def}{Definition}[section]
\newtheorem{ex}[Def]{Example}
\newtheorem{rem}[Def]{Remark}

\theoremstyle{plain}
\newtheorem{prop}[Def]{Proposition}
\newtheorem{thm}[Def]{Theorem}
\newtheorem*{thm*}{Theorem}
\newtheorem{lem}[Def]{Lemma}
\newtheorem{cor}[Def]{Corollary}
\newtheorem*{cor*}{Corollary}
\newtheorem{con}[Def]{Conjecture}
\newtheorem*{con*}{Conjecture}

\newtheorem*{verm*}{Vermutung}

\usepackage[draft]{fixme}
\fxsetup{theme=color, mode=multiuser, noinline}
\FXRegisterAuthor{mk}{MK}{\color{purple}Mario}
\FXRegisterAuthor{ds}{DS}{\color{teal}David}

\include{operators}
\include{fonts}

\title[Three results on the half-plane property]{Three results related to the half-plane property of matroids}
\author{Mario Kummer}
\address{Technische Universit\"at, Dresden, Germany} 
\email{mario.kummer@tu-dresden.de}
\author{David Sawall}
\address{Universit\"at Konstanz}
\email{david.sawall@uni-konstanz.de}
\thanks{Mario Kummer has been supported by the Deutsche Forschungsgemeinschaft under Grant No. 502861109.}
\begin{document}


\begin{abstract}
 We settle three problems from the literature on stable and real zero polynomials and their connection to matroid theory. We disprove the weak real zero amalgamation conjecture by Schweighofer and the second author. We disprove a conjecture by Br\"and\'en and D'Le\'on by finding a relaxation of a matroid with the weak half-plane property that does not have the weak half-plane property itself. Finally, we prove that every quaternionic unimodular matroid has the half-plane property which was conjectured by Pendavingh and van Zwam.
\end{abstract}
\maketitle

\section{Introduction}
A \emph{real zero polynomial} is a polynomial $P\in\R[x_1,\ldots,x_n]$ with $P(0)\neq0$ which has only real zeros when restricted to any real line through the origin. Real zero polynomials play an important role in the theory of semidefinite programming, see e.g. \cite{victorsurvey}. In this context, the second author and Schweighofer made the following conjecture:

\begin{con*}[{\cite[Conjecture~7.6]{rzamalgamation}}]
    Assume that $F\in\R[x_1,x_2,y_1,\ldots,y_m]$ and $G\in\R[x_1,x_2,z_1,\ldots,z_n]$ are real zero polynomials. If
    \begin{equation*}
        F|_{y_1=\cdots=y_m=0}=G|_{z_1=\cdots=z_n=0},
    \end{equation*}
    then there is a real zero polynomial $H\in\R[x_1,x_2,y_1,\ldots,y_m,z_1,\ldots,z_n]$ such that
    \begin{equation*}
        F=H|_{z_1=\cdots=z_n=0}\textrm{ and }G=H|_{y_1=\cdots=y_m=0}.
    \end{equation*}
\end{con*}

We provide a counterexample to this conjecture in \Cref{sec:amal}. Using the connection between stable polynomials, a concept closely related to real zero polynomials, and the theory of matroids and polymatroids established in \cite{HPPlong} and \cite{BrandenHPP}, we associate a polymatroid to every real zero polynomial. Then the key step is \Cref{thm:amalmain} which says that if two real zero polynomials satisfy the conclusion the above conjecture, then the associated polymatroids can be \emph{amalgamated} in the sense of \Cref{def:amal}. Our counterexample comes from two polymatroids which cannot be amalgamated. These polymatroids are obtained by specializing the matroids $F_7^{-4}$ and $F_7^{-5}$ --- two relaxations of the Fano matroid which have the \emph{half-plane property}, meaning that their bases generating polynomials are stable. 

More generally, a matroid has the \emph{weak half-plane property} if its set of bases agrees with the support of a stable polynomial. The question of which matroids have the (weak) half-plane property has been extensively studied. In this context, Br\"and\'en and D'Le\'on offered the following conjecture:

\begin{con*}[{\cite[Conjecture~4.2]{BrandenDLeon}}]
    Suppose that $M$ has the weak half-plane property. Then so does any relaxation of $M$.
\end{con*}

In \Cref{sec:whpp} we present a counterexample to this conjecture. More precisely, we show that a relaxation $P_1$ of the matroid $P_8$ does not have the weak half-plane property. On the other hand, $P_8$ does have the weak half-plane property because it is representable over $\R$. For proving that $P_1$ does not have the weak half-plane property, we employ the techniques developed in \cite{BrandenDLeon} for narrowing down the space of possible coefficients of a hypothetical stable polynomial with support $P_1$.

A key feature of the matroid $P_1$ that made us examine it for the weak half-plane property is that it is not representable --- even under more general notions of representability as discussed in \cite{skewpartial}. Matroids representable over more general structures than fields still often tend to have the weak half-plane property, see e.g. \cite{NonPappus}. Our third result is of this flavor. Namely, we prove in \Cref{sec:qu} a conjecture which was attributed by Pendavingh and van Zwam to David G. Wagner:

\begin{con*}[{\cite[Conjecture~6.9]{skewpartial}}]
    All quaternionic unimodular matroids have the half-plane property.
\end{con*}

The notion of quaternionic unimodular matroids is a generalization of the class of sixth root of unity matroids to the skew field of quaternions. For a precise defintion see \Cref{def:qu}.

\section{Preliminaries}
We denote by $\N$ and $\N_0$ the set of positive and nonnegative integers, respectively.
In this section we let $E$ always denote a finite set. For $i\in E$ we denote by $\delta_i\in\R^E$ the $i$th unit vector. For $x\in\R^E$ we write $|x|=\sum_{i\in E}|x_i|$. We further let $[n]=\{1,\ldots,n\}$ and denote by $\binom{[n]}{k}$ the set of all $k$-element subsets of $[n]$ for every $k,n\in\N$. We recall the cryptomorphic definitions of M-convex sets and polymatroids.

\begin{Def}
     A subset $J\subseteq\N_0^E$ is \emph{M-convex} if for every $i\in E$ and every $\alpha,\beta\in J$ such that $\alpha_i>\beta_i$, there is $j\in E$ satisfying
            \begin{equation*}
                \alpha_j<\beta_j\textrm{ and }\alpha-\delta_i+\delta_j\in J.
            \end{equation*}
\end{Def}

\begin{Def}
    A \emph{polymatroid} on $E$ is a function $r\colon 2^E\to\N_0$ such that we have for all $S,T\subseteq E$:
            \begin{enumerate}[(i)]
            \item $r(\emptyset)=0$,
                \item $r(S)\leq r(T)$ if $S\subseteq T$, and
                \item $r(S\cup T)+r(S\cap T)\leq r(S)+r(T)$.
            \end{enumerate}
\end{Def}
If $J\subseteq\N_0^E$ is an M-convex set, then we define the function $r_J\colon2^E\to\N_0$ by
    \begin{equation*}
        r_J(S)=\max\{\sum_{i\in S}\alpha_i\mid \alpha\in J\}.
    \end{equation*}
Conversely, if $r\colon 2^E\to\N_0$ is a polymatroid, then we define the set $J_r\subseteq\N_0^E$ by
\begin{equation*}
    J_r=\{x\in\N_0^E\mid \sum_{i\in S}x_i\leq r(S)\textrm{ for all }S\subseteq E\textrm{ and }\sum_{i\in E}x_i=r(E)\}.
\end{equation*}
These two constructions are inverse to each other and define bijections between the set of M-convex sets in $\N_0^E$ and the set of polymatroids on $E$, see for example \cite[\S4.4]{murota}.

\begin{rem}
    A polymatroid $r\colon 2^E\to\N_0$ is the rank function of a matroid $M$ on $E$ if and only if $r(\{i\})\leq1$ for all $i\in E$. In this case we have
    \begin{equation*}
        J_r=\{\sum_{i\in B}\delta_i\mid B\textrm{ is a basis of }M\}.
    \end{equation*}
    We apply the definitions made for polymatroids to matroids by considering their rank functions as polymatroids.
\end{rem}

\subsection{Amalgamation of polymatroids}
There has been considerable interest in the question whether two (poly)matroids can be amalgamated in the following sense, see for example \cite{sticky,stickypolymatroids} and also \cite[\S11.4]{oxley}.
\begin{Def}\label{def:amal}
    Let $E_1,E_2$ be two finite sets, let $E_0=E_1\cap E_2$ and $E_3=E_1\cup E_2$. Let $r_i$ be a polymatroid on $E_i$ for all $i=0,1,2,3$.
    \begin{enumerate}[(a)]
        \item We say that $r_0$ is the \emph{restriction} of $r_1$ to $E_0$ if $r_0(S)=r_1(S)$ for all $S\subseteq E_0$ and denote this by $r_0=r_1|_{E_0}$.
        \item We say that $r_3$ is an \emph{amalgam} of $r_1$ and $r_2$ if $r_1=r_3|_{E_1}$ and $r_2=r_3|_{E_2}$.
    \end{enumerate}
\end{Def}

Clearly, a necessary condition for an amalgam of $r_1$ and $r_2$ as in \Cref{def:amal} to exist is that $r_1|_{E_0}=r_2|_{E_0}$. However, this condition is not sufficient \cite[Example~7.2.4]{oxley}.

\begin{Def}
    A polymatroid $r_0$ on $E$ is called \emph{sticky} if there exists an amalgam for all polymatroids $r_1$ and $r_2$ on finite sets $E_1$ and $E_2$ with $E=E_1\cap E_2$ and such that $r_0=r_i|_{E}$ for $i=1,2.$
\end{Def}

In order to describe conditions for a polymatroid to be sticky, recall the following definition, generalizing the usual notion for matroids.
\begin{Def}
    Let $r$ be a polymatroid on $E$.
    \begin{enumerate}[(a)]
        \item A subset $F\subseteq E$ is called a \emph{flat} of $r$ if $r(F')>r(F)$ for every proper superset $F'$ of $F$.
        \item Two flats $F_1,F_2$ are called a \emph{modular pair} if
        \begin{equation*}
            r(F_1)+r(F_2)=r(F_1\cup F_2)+r(F_1\cap F_2).
        \end{equation*}
    \end{enumerate}
\end{Def}

\begin{thm}[\cite{stickypolymatroids}]\label{thm:stickypoly}
    Let $r$ be a polymatroid on $E$. If every pair of flats is modular, then $r$ is sticky. The converse holds if $|E|\leq5$.
\end{thm}

\begin{ex}\label{ex:uniform}
 Let $2\leq n\leq5$ and $E=[n]$. One checks that
 \begin{equation*}
     J=\{x\in\N_0^n\mid |x|=3\textrm{ and }\forall i\in E: x_i\leq 2\}
 \end{equation*}
 is M-convex. The polymatroid $r_J$ is not sticky because the pair $\{1\},\{2\}$ of flats is not modular.
\end{ex}

\subsection{Stable polynomials}
We briefly recall the different stability and real zero properties of polynomials and their relation to each other. As a general reference we recommend the survey \cite{pemantle}.
\begin{Def}
    Let $P\in\R[x_i\mid i\in E]$ be a polynomial.
    \begin{enumerate}[(a)]
        \item The \emph{support} of $P$ is the unique subset $\supp(P)\subseteq\N_0^E$ such that we can write
        \begin{equation*}
            P=\sum_{\alpha\in\supp(h)}c_\alpha x^\alpha
        \end{equation*}
        for some non-zero $c_\alpha\in\R$.
        \item The polynomial $P$ is called \emph{stable} if for all $z\in\C^E$ such that $\Ima(z_i)>0$ for all $i\in E$ we have $P(z)\neq0$.
        \item The polynomial $P$ is called a \emph{real zero polynomial} if for all $v\in\R^E$ the univariate polynomial $P(t\cdot v)\in\R[t]$ has only real zeros.
        \item Let $P$ be homogeneous and $e\in\R^E$. Then $P$ is called \emph{hyperbolic} with respect to $e$ if for all $v\in\R^E$ the univariate polynomial $P(t\cdot e+v)\in\R[t]$ has only real zeros.
    \end{enumerate}
\end{Def}

Stability is preserved under taking partial derivatives in coordinate directions, under setting some variables equal to each other and under scaling variables by positive scalars. We summarize the well-known connection between the above concepts.

\begin{prop}[see for example {\cite[Proposition 5.3]{pemantle}}]\label{prop:stablehyp}
    Let $P\in\R[x_i\mid i\in E]$ be a homogeneous polynomial. The following are equivalent:
    \begin{enumerate}[(i)]
        \item $P$ is stable.
        \item $P$ is hyperbolic with respect to every point in $\R_{>0}^E$.
        \item $P$ is hyperbolic with respect to every point in $\R_{>0}^E$ and every $e\in\R_{\geq0}^E$ with $P(e)\neq0$.
    \end{enumerate}
\end{prop}

\begin{lem}
    Let $P\in\R[x_i\mid i\in E]$ be a homogeneous polynomial and $i\in E$. The following are equivalent:
    \begin{enumerate}[(i)]
        \item $P$ is hyperbolic with respect to $\delta_i$.
        \item $P|_{x_i=1}\in\R[x_i\mid i\in E\smallsetminus\{i\}]$ is a real zero polynomial.
    \end{enumerate}
\end{lem}

\begin{Def}
    Let $P\in\R[x_i\mid i\in E]$ be hyperbolic with respect to $e\in\R^E$. The \emph{hyperbolicity cone} of $P$ at $e$ is defined as
    \begin{equation*}
        \Lambda(P,e)=\{v\in\R^E\mid P(t\cdot e+v)\in\R[t]\textrm{ has only nonnegative zeros}\}.
    \end{equation*}
\end{Def}

Hyperbolicity cones are convex cones \cite{gar}. \Cref{prop:stablehyp} can be rephrased in terms of hypberbolicity cones.

\begin{prop}
   Let $P\in\R[x_i\mid i\in E]$ be a homogeneous polynomial and $e\in\R_{>0}^E$. The following are equivalent:
   \begin{enumerate}[(i)]
        \item $P$ is stable.
       \item $P$ is hyperbolic with respect to $e$ and $\delta_i\in\Lambda(P,e)$ for all $i\in E$.
   \end{enumerate}
\end{prop}

The connection of stable polynomials to M-convex sets is given by the following.

\begin{thm}[{\cite[Theorem~3.2]{BrandenHPP}}]
    Let $P\in\R[x_i\mid i\in E]$ be a homogeneous stable polynomial. Then $\supp(P)$ is M-convex.
\end{thm}

\begin{Def}
    Let $J\subseteq\N_0^E$ be an M-convex set.
    \begin{enumerate}[(a)]
        \item The \emph{generating polynomial} of $J$ is
        \begin{equation*}
            h_J=\sum_{\alpha\in J}\frac{x^\alpha}{\alpha!}\in \R[x_i\mid i\in E]
        \end{equation*}
        where $\alpha!=\prod_{i\in E}\alpha_i!$.
        \item We say that $J$ has the \emph{half-plane property} if its generating polynomial is stable.
        \item We say that $J$ has the \emph{weak half-plane property} if there exists a homogeneous stable polynomial $P\in\R[x_i\mid i\in E]$ such that $J=\supp(P)$.
        \item We say that a polymatroid $r$ on $E$ has the \emph{(weak) half-plane property} if the associated M-convex set $J_r$ has the (weak) half-plane property.
    \end{enumerate}
\end{Def}

In the remaining part, we observe that in certain nice situations, restricting the polymatroid of a homogeneous stable polynomial corresponds to plugging in zeros for some of the variables.

\begin{Def}
    Let $J\subseteq\N_0^E$ be M-convex and let $T\subseteq E$. We say that $J$ is \emph{nondegenerate with respect to $T$} if there exists $\alpha\in J$ such that for all $i\in E\smallsetminus T$ we have $\alpha_i=0$.
\end{Def}

\begin{lem}\label{lem:mconvrest}
    Let $J\subseteq\N_0^E$ be M-convex and nondegenerate with respect to $T\subseteq E$. Then we have $r_J|_T=r_{J'}$ where
    \begin{equation*}
        J'=\{(\alpha_i)_{i\in T}\mid \alpha\in J\textrm{ such that }\forall i\in E\smallsetminus T: \alpha_i=0\}\subseteq\N_0^T.
    \end{equation*}
\end{lem}

\begin{proof}
    Let $S\subseteq T$. Then
    \begin{equation*}
        r_J|_T(S)=r_J(S)=\max\{\sum_{k\in S}\alpha_k\mid \alpha\in J\}.
    \end{equation*}
    Let $\alpha\in J$ be a point where this maximum is attained such that $\sum_{k\in E\smallsetminus T}\alpha_k$ is minimal. If $\sum_{k\in E\smallsetminus T}\alpha_k=0$, then $\alpha\in J'$ and we have $r_J|_T(S)=r_{J'}(S)$. Thus assume for the sake of a contradiction that $\alpha_i>0$ for some $i\in E\smallsetminus T$. Since $J$ is nondegenerate with respect to $T\subseteq E$, there exists $\beta\in J$ such that for all $k\in E\smallsetminus T$ we have $\beta_k=0$. Because $J$ is M-convex, there exists $j\in T$ such that $\alpha_j<\beta_j$ and $\gamma=\alpha-\delta_i+\delta_j\in J$. Then $\sum_{k\in S}\gamma_k\geq\sum_{k\in S}\alpha_k$ because $i\not\in S$ and $\sum_{k\in E\smallsetminus T}\gamma_k<\sum_{k\in E\smallsetminus T}\alpha_k$ because $i\in E\smallsetminus T$ and $j\in T$. This contradicts our minimality assumption.
\end{proof}

\begin{cor}\label{cor:polyamal}
    Let $E=E_1\cup E_2$ and let $P\in\R[x_i\mid i\in E]$ be a polynomial such that $J=\supp(P)$ is M-convex.
    Finally, for $k=1,2$, we let
    \begin{equation*}
        P_k=P|_{x_i=0\textrm{ for }i\not\in E_k}\in \R[x_i\mid i\in E_k]
    \end{equation*}
    and $J_k=\supp(P_k)$. If $P_1,P_2$ are both not identically zero, then $J_1$ and $J_2$ are M-convex and $r_J$ is an amalgam of $r_{J_1}$ and $r_{J_2}$.
\end{cor}

\begin{proof}
    The condition that $P_k$ is not identically zero is equivalent to $J$ being nondegenerate with respect to $E_k$. Now the claim follows from \Cref{lem:mconvrest}.
\end{proof}

\section{Real zero amalgamation}\label{sec:amal}
In this section we always let $E=E_1\cup E_2$ be a finite set and $E_0=E_1\cap E_2$. We assume that $0\in E_0$. For $k=0,1,2$ we let $0\neq P_k\in\R[x_i\mid i\in E_k]$ be homogeneous and stable such that
\begin{equation*}
    P_0=P_k|_{x_i=0\textrm{ for }i\not\in E_0}
\end{equation*}
for $k=1,2$. Note that this implies that $P_0,P_1$ and $P_2$ all have the same degree $d$. We further denote $J_k=\supp(P_k)$ and $r_k=r_{J_k}$ for $k=0,1,2$. For $k=0,1,2$ we also consider the polynomial $H_k$ obtained from $P_k$ by substituting $x_0+x_i$ for $x_i$ for all $i\in E_0\smallsetminus\{0\}$.

\begin{lem}\label{lem:hypcone1}
    The polynomial $H_k$ is hyperbolic with respect to $\delta_0$. Its hyperbolicity cone contains $\delta_i$ for all $i\in E_k$ and the point
    \begin{equation*}
        \delta_0-\sum_{0\neq i\in E_0}\delta_i.
    \end{equation*}
\end{lem}

\begin{proof}
    The statement of $H_k$ being hyperbolic with respect to $\delta_0$ is equivalent to $P_k$ being hyperbolic with respect to $\sum_{i\in E_0}\delta_i$. Because of
    \begin{equation*}
        P_k\left(\sum_{i\in E_0}\delta_i\right)=P_0\left(\sum_{i\in E_0}\delta_i\right)\neq0,
    \end{equation*}
    by stability of $P_0$, this follows because $P_k$ is stable. The hyperbolicity cone of $H_k$ containing $\delta_i$, $i\in E_k$, follows from the same statement for $P_k$. The last statement is equivalent to $\delta_0$ being in the hyperbolicity cone of $P_k$.
\end{proof}

\begin{lem}\label{lem:hypcone2}
    Let $H\in\R[x_i\mid i\in E]$ be hyperbolic with respect to $\delta_0$ such that the hyperbolicity cone of
    \begin{equation*}
     H|_{x_i=0\textrm{ for }i\not\in E_k}
    \end{equation*}
    agrees with the hyperbolicity cone of $H_k$ for $k=1,2$. Then the hyperbolicity cone of $H$ contains $\delta_i$ for all $i\in E$ and the point
    \begin{equation*}
        \delta_0-\sum_{0\neq i\in E_0}\delta_i.
    \end{equation*}
\end{lem}

\begin{proof}
    This follows immediately from \Cref{lem:hypcone1}.
\end{proof}

Now let $E'=E\smallsetminus\{0\}$ and likewise $E_k'=E_k\smallsetminus\{0\}$ for $k=0,1,2$. Then we consider the polynomials
\begin{equation*}
    Q_k=H_k|_{x_0=1}\in\R[x_i\mid i\in E_k']
\end{equation*}
for $k=0,1,2$. These are real zero polynomials with the property that
\begin{equation*}
    Q_0=Q_k|_{x_i=0\textrm{ for }i\not\in E'_0}
\end{equation*}
or $k=1,2$. We have the following result.

\begin{thm}\label{thm:amalmain}
    Assume that there is a real zero polynomial $Q\in\R[x_i\mid i\in E']$ with 
    \begin{equation*}
        Q_k=Q|_{x_i=0\textrm{ for }i\not\in E'_k}
    \end{equation*}
    for $k=1,2$. Then $r_1$ and $r_2$ have an amalgam.
\end{thm}

\begin{proof}
    Let $d'=\deg(Q)$. We define
    \begin{equation*}
        H=x_0^{d'}\cdot Q(\frac{x_i}{x_0}\mid i\in E')\in\R[x_i\mid i\in E].
    \end{equation*}
    Then we have
    \begin{equation*}
        H|_{x_i=0\textrm{ for }i\not\in E'_k}=x_0^{d'-d}\cdot H_k
    \end{equation*}
    for $k=1,2$. Thus by \Cref{lem:hypcone2} the hyperbolicity cone of $H$ contains $\delta_i$, $i\in E$, and the point
    \begin{equation*}
        \delta_0-\sum_{0\neq i\in E_0}\delta_i.
    \end{equation*}
    This implies that the polynomial $P$ obtained from $H$ by substituting $x_i-x_0$ for $x_i$ for all $i\in E_0\smallsetminus\{0\}$ is stable. We have
    \begin{equation*}
        P|_{x_i=0\textrm{ for }i\not\in E_k}=x_0^{d'-d}\cdot P_k
    \end{equation*}
    for $k=1,2$. Finally, we let $P'$ be the polynomial obtained from $P$ by dropping all monomials that are not divisible by $x_0^{d'-d}$ and dividing the result by $x_0^{d'-d}$. Then
    \begin{equation*}
        P'|_{x_i=0\textrm{ for }i\not\in E_k}=P_k
    \end{equation*}
    for $k=1,2$. The support of $P'$ is M-convex because it agrees with the support of the stable polynomial $\frac{\partial^{d'-d}}{\partial x_0}P$. Now the claim follows from \Cref{cor:polyamal}.
\end{proof}

Now we are ready to disprove the real zero amalgamation conjecture from \cite{rzamalgamation}.

\begin{con}[{\cite[Conjecture~7.6]{rzamalgamation}}]\label{con:weakrzamal}
    Let $E'=E_1'\cup E_2'$ be a finite set such that $E_0'=E_1'\cap E_2'$ has two elements. For $k=1,2$ let $Q_k\in\R[x_i\mid i\in E_k']$ be a real zero polynomial.
    If
    \begin{equation*}
        Q_1|_{x_i=0\textrm{ for }i\not\in E_0'}=Q_2|_{x_i=0\textrm{ for }i\not\in E_0'},
    \end{equation*}
    then there is a real zero polynomial $Q\in\R[x_i\mid i\in E']$ such that 
    \begin{equation*}
     Q_k=Q|_{x_i=0\textrm{ for }i\not\in E_k'}   
    \end{equation*}
     for $k=1,2$.
\end{con}

\Cref{thm:amalmain} shows that for $|E_0|=3$ (and hence $|E_0'|=2$) every choice of $P_k$, $k=0,1,2$, as above such that additionally $r_1$ and $r_2$ do not have an amalgam gives a counterexample to \Cref{con:weakrzamal}.

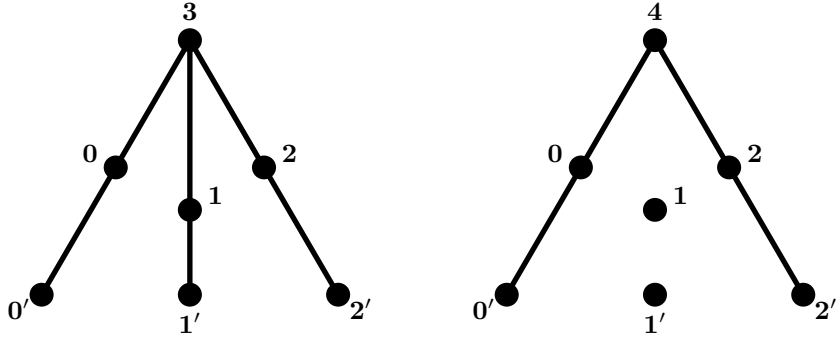
\begin{figure}
        \begin{tikzpicture}[line cap=round,line join=round,>=triangle 45,x=.75cm,y=0.75cm]
\clip(-4,-4) rectangle (4,4);
\draw [line width=2pt] (0,3)-- (-2.6,-1.5);
\draw [line width=2pt] (0,3)-- (2.6,-1.5);
\draw [line width=2pt] (0,3)-- (0,-1.5);
\fill [color=black] (0,3) circle (4.5pt);
\draw[color=black] (0,3.5) node {\textbf{3}};
\fill [color=black] (0,0) circle (4.5pt);
\draw[color=black] (0.45,0.25) node {\textbf{1}};
\fill [color=black] (-2.6,-1.5) circle (4.5pt);
\draw[color=black] (-3,-1.75) node {\textbf{0$'$}};
\fill [color=black] (2.6,-1.5) circle (4.5pt);
\draw[color=black] (3,-1.75) node {\textbf{2$'$}};
\fill [color=black] (-1.3,0.75) circle (4.5pt);
\draw[color=black] (-1.75,1) node {\textbf{0}};
\fill [color=black] (1.3,0.75) circle (4.5pt);
\draw[color=black] (1.75,1) node {\textbf{2}};
\fill [color=black] (0,-1.5) circle (4.5pt);
\draw[color=black] (0,-2) node {\textbf{1$'$}};
\end{tikzpicture}
\begin{tikzpicture}[line cap=round,line join=round,>=triangle 45,x=0.75cm,y=0.75cm]
\clip(-4,-4) rectangle (4,4);
\draw [line width=2pt] (0,3)-- (-2.6,-1.5);
\draw [line width=2pt] (0,3)-- (2.6,-1.5);
\fill [color=black] (0,3) circle (4.5pt);
\draw[color=black] (0,3.5) node {\textbf{4}};
\fill [color=black] (0,0) circle (4.5pt);
\draw[color=black] (0.45,0.25) node {\textbf{1}};
\fill [color=black] (-2.6,-1.5) circle (4.5pt);
\draw[color=black] (-3,-1.75) node {\textbf{0$'$}};
\fill [color=black] (2.6,-1.5) circle (4.5pt);
\draw[color=black] (3,-1.75) node {\textbf{2$'$}};
\fill [color=black] (-1.3,0.75) circle (4.5pt);
\draw[color=black] (-1.75,1) node {\textbf{0}};
\fill [color=black] (1.3,0.75) circle (4.5pt);
\draw[color=black] (1.75,1) node {\textbf{2}};
\fill [color=black] (0,-1.5) circle (4.5pt);
\draw[color=black] (0,-2) node {\textbf{1$'$}};
\end{tikzpicture}
        \caption{The matroids $F_7^{-4}$ (left) and $F_7^{-5}$ (right).}
        \label{fig:f7relax}
 \end{figure}

\subsection{An example with $|E_0|=3$}
We consider the matroids $F_7^{-4}$ and $F_7^{-5}$ defined as in \Cref{fig:f7relax}. These matroids both have the half-plane property, see \cite[\S A.2.2]{HPPlong} for $F_7^{-5}$ and \cite{wagnerwei} for $F_7^{-4}$. This means that their generating polynomials $h_{F_7^{-4}}$ and $h_{F_7^{-5}}$ are stable. We let $P_1$ and $P_2$ be the polynomials obtained from $h_{F_7^{-4}}$ and $h_{F_7^{-5}}$, respectively, by setting $x_{0'}=x_0$, $x_{1'}=x_1$ and $x_{2'}=x_2$. By construction we have
\begin{equation*}
    P_1|_{x_3=0}=P_2|_{x_4=0}.
\end{equation*}
Hence $E_1=\{0,1,2,3\}$ and $E_2=\{0,1,2,4\}$. Note that the support of $P_0:=P_1|_{x_3=0}$ is the M-convex set from \Cref{ex:uniform} for $n=3$. We will show that the polymatroids $r_1$ and $r_2$ corresponding to $P_1$ and $P_2$ do not have an amalgam so that the real zero polynomials $Q_1$ and $Q_2$ obtained from $P_1$ and $P_2$ constitute a counterexample to \Cref{con:weakrzamal}. In fact, we will even prove that for every $m\in\N$ the polymatroids $m\cdot r_1$ and $m\cdot r_2$ do not have an amalgam which shows that $Q_1^m$ and $Q_2^m$ also do not satisfy the conclusion of \Cref{con:weakrzamal} for any $m\in\N$.

\begin{thm}
    The polymatroids $m\cdot r_1$ and $m\cdot r_2$ do not have an amalgam.
\end{thm}

\begin{proof}
    We proceed as in the proof of \cite[Theorem~2]{sticky}.
    Assume for the sake of a contradiction that the polymatroid $r$ on $E=\{0,1,2,3,4\}$ is an amalgam of $m\cdot r_1$ and $m\cdot r_2$. We have
    \begin{equation*}
        r(\{0\})+r(\{0,3\})\leq r(\{0\})+r(\{0,3,4\})\leq r(\{0,3\})+r(\{0,4\}).
    \end{equation*}
    By definition of $r_1$ and $r_2$ this implies that
    \begin{equation*}
        2m+2m\leq 2m+r(\{0,3,4\})\leq 2m+2m.
    \end{equation*}
    Thus we have $r(\{0,3,4\})=2m$. Likewise one shows that $r(\{2,3,4\})=2m$. Furthermore, we have that
    \begin{equation*}
      4m\leq r(\{3,4\})+3m \leq r(\{3,4\})+r(\{0,2,3,4\})\leq r(\{0,3,4\})+r(\{2,3,4\})=4m.
    \end{equation*}
    This shows $r(\{3,4\})=m$. Finally, we have
    \begin{equation*}
      m+2m\leq  r(\{3\})+r(\{1,3,4\}) \leq r(\{1,3\})+r(\{3,4\}) =2m+m
    \end{equation*}
    which implies $r(\{1,3,4\})=2m$ contradicting
    \begin{equation*}
        3m=r(\{1,4\})\leq r(\{1,3,4\}).
     \end{equation*}
     Hence the polymatroids $m\cdot r_1$ and $m\cdot r_2$ do not have an amalagam.
\end{proof}

\section{Weak half-plane property}\label{sec:whpp}
Let $E$ always denote a finite set and $M$ a matroid on $E$ with set of bases $\cB$.

\begin{Def}[\cite{BrandenDLeon}]
    We say that four bases $B_1,B_2,B_3,B_4\in\cB$ form a \emph{degenerate quadrangle} of $M$ if there exists $S\subseteq E$ and pairwise different $i,j,k,l\notin S$ such that 
    \begin{equation*}
      (B_1,B_2,B_3,B_4)=(S\cup\{i,k\},S\cup\{j,l\},S\cup\{i,l\},S\cup\{j,k\})  
    \end{equation*}
    and if at most one of $S\cup\{i,j\}$ and $S\cup\{k,l\}$ is a bases of $M$. 
\end{Def}

The following theorem was used in \cite{BrandenDLeon} to reduce the number of possible parameters when searching for a stable polynomial with support $M$.

\begin{thm}[\cite{BrandenHPP}]
    For every basis $B\in\cB$ let $0\neq a_B\in\R$ be such that the multiaffine and homogeneous polynomial
    \begin{equation*}
        P=\sum_{B\in\cB}a_{B}\cdot\prod_{i\in B}x_i\in\R[x_i\mid i\in E]
    \end{equation*}
    is stable. If $B_1,B_2,B_3,B_4$ form a degenerate quadrangle of $M$, then
    \begin{equation}\label{eq:deg4}
      a_{B_1}a_{B_3}=a_{B_2}a_{B_4}.  
    \end{equation}
\end{thm}
 Letting $b_B:=\log(|a_B|)$ for all $B\in\cB$ we obtain from \Cref{eq:deg4} linear equations
 \begin{equation*}
  b_{B_1}+b_{B_3}-b_{B_2}-b_{B_4}=0   
 \end{equation*}
for all degenerate quadrangles $B_1,B_2,B_3,B_4$ of $M$. We denote by $V_M\subseteq\R^{\cB}$ the linear space cut out by all such equations. By \cite[Lemma~2.6]{BrandenDLeon} the vector space
\begin{equation*}
  W_M:=\{(\sum_{i\in B}v_i)_{B\in\cB}\mid v\in\R^E\}  
\end{equation*}
is an $(|E|-z+1)$-dimensional linear subspace of $V_M$ where $z$ is the number of connected components of $M$.

\begin{lem}\label{lem:rescale}
    Let $U_M$ be a linear complement of $W_M$ in $V_M$. If $M$ has the weak half-plane property, then there exists a vector $b\in U_M$ such that
    \begin{equation*}
        \sum_{B\in\cB}\exp(b_{B})\cdot\prod_{i\in B}x_i
    \end{equation*}
    is stable.
\end{lem}

\begin{proof}
    This follows in the same way as \cite[Theorem~2.3]{BrandenDLeon}: Scaling the variables by $x_i\mapsto\exp(v_i)x_i$ corresponds, after taking logarithms of the coefficients, to shifting by the corresponding vector from $W_M$.
\end{proof}

For choosing a linear complement of $W_M$ in $V_M$ in a nice way, the following lemma might be useful.

\begin{lem}\label{lem:notreg}
    Let $M$ be represented by a matrix $A\in\R^{d\times |E|}$ of rank $d$. For $B\in\cB$ we denote by $A[B]$ the corresponding $d\times d$ submatrix. Then 
    \begin{equation*}
        u(A):=(\log|\det(A[B])|)_{B\in\cB}\in V_M.
    \end{equation*}
    Furthermore, if $M$ is not regular, then $u(A)\not\in W_M$.
\end{lem}

\begin{proof}
    By \cite[Theorem~8.1]{HPPlong} the polynomial
    \begin{equation*}
        \sum_{B\in\cB}\det(A[B])^2\cdot\prod_{i\in B}x_i
    \end{equation*}
    is stable. This proves the first claim. Now assume that there exists $v\in\R^E$ such that $u(A)=(\sum_{i\in B}v_i)_{B\in\cB}$. Scaling the $i$th column of $A$ by $\exp(-v_i)$ for all $i\in E$, we obtain a matrix $A'$ representing $M$ all of whose maximal minors are in $\{-1,0,1\}$. After multiplication of $A'$ from the left by a suitable invertible matrix we can additionally assume that $A'[B_0]$ is the identity matrix for some $B_0\in\cB$. Then $A'$ is a totally unimodular matrix representing $M$ which shows that $M$ is regular.
\end{proof}

\begin{ex}\label{ex:p8}
    In this example, we consider the rank 4 matroid $M=P_8$ on 8 elements which is represented by the real matrix
    \[A:=\begin{pmatrix}
        1&0&0&0&0&1&1&2\\
        0&1&0&0&1&0&1&1\\
        0&0&1&0&1&1&0&1\\
        0&0&0&1&2&1&1&0
    \end{pmatrix}\]
    whose columns we label by $0,\ldots,7$.
    This matroid is not regular \cite[\S A4]{HPPlong} and therefore $u(A)$ is in $V_M$ but not in $W_M$ by \Cref{lem:notreg}. Using the {\texttt{Macaulay2}} \cite{M2} package ``Matroids'' \cite{MatroidsArticle} we compute that $\dim(V_M)=9$. Because $M$ is connected, this implies that the span of $u(A)$ is a linear complement of $W_M$ in $V_M$. Note that $P_8$ has the weak half-plane property because it is representable over $\R$. 
\end{ex}

Recall that if $X$ is a circuit-hyperplane of $M$, then $\cB\cup\{X\}$ is the set of bases of a matroid $M'$ \cite[Theorem~1.5.14]{oxley}. Then $M'$ is called a \emph{relaxation} of $M$.

\begin{lem}\label{lem:relax}
    Let $M'$ be a relaxation of $M$. Then every degenerate quadrangle of $M'$ is a degenerate quadrangle of $M$.
\end{lem}

\begin{proof}
    We proceed as in the proof of \cite[Lemma~3.43]{masterleon}.
    Denote by $X$ the circuit-hyperplane of $M$ such that $M'$ is the relaxation of $M$ by $X$.
    We show that for all $x\in X$ and $y\in E\smallsetminus X$ the set
    \[(X\smallsetminus\{x\})\cup\{y\}\]
    forms a basis of $M$ (and therefore of $M'$). Then it immediately follows that $X$ cannot be contained in some degenerate quadrangle of $M'$.
    Since $X$ is a circuit, $X\smallsetminus\{x\}$ still has rank $\rank(M)-1=\rank(M')-1$. Then $(X\smallsetminus\{x\})\cup\{y\}$ has rank $\rank(M')$ because $X$ is closed and $y\notin X$. Thus $(X\smallsetminus\{x\})\cup\{y\}$ is a basis of $M'$.
\end{proof}

\begin{rem}\label{rem:embed}
    Assume that $\rank(M)\geq2$.
    Let $M'$ be the relaxation of $M$ by the circuit-hyperplane $X$ and denote $\cB'=\cB\cup\{X\}$. \Cref{lem:relax} implies that the map $\R^{\cB}\to\R^{\cB'}$, that sends $v\in\R^{\cB}$ to the vector $w$ with $w_{X}=0$ and $w_B=v_B$ for all $B\in\cB$, maps $V_M$ to $V_{M'}$. We thus obtain a natural embedding
    \begin{equation*}
        \iota\colon V_M\hookrightarrow V_{M'}.
    \end{equation*}
    Furthermore, we have $\delta_X\in V_{M'}$. If $U_M$ is a linear complement of $W_M$ in $V_M$, then
    \begin{equation*}
        \left(\iota(U_M)\oplus\R\cdot\delta_X\right)\cap W_{M'}=\{0\}.
    \end{equation*}
    However, in general $\iota(U_M)$ and $\delta_X$ do not span a linear complement of $W_{M'}$, see for example \cite[Table~1]{BrandenDLeon}.
\end{rem}

\begin{ex}\label{ex:f7}
    As an illustration, we recall \cite[Example~4.1]{BrandenDLeon}. The non-Fano matroid $M'=F_7^{-}$ is represented by the real matrix
    \[A=\begin{pmatrix}
        1&1&0&0&0&1&1\\
        0&1&1&1&0&0&1\\
        0&0&0&1&1&1&1
    \end{pmatrix}\]
    whose columns we label by $1,\ldots,7$. This matroid is not regular \cite[\S A.2.2]{HPPlong} and therefore $u(A)$ is in $V_{M'}$ but not in $W_{M'}$ by \Cref{lem:notreg}. On the other hand, the non-Fano matroid is a relaxation of the Fano matroid $F_7$ by the circuit hyperplane $\{2,4,6\}$. Because $\dim(V_{M'}/W_{M'})=1$ by \cite{BrandenDLeon}, it follows from \Cref{rem:embed} that $u(A)+W_{M'}$ must contain a scalar multiple of $\delta_{\{2,4,6\}}$. Indeed, one has that $u(A)=\log(2)\cdot\delta_{\{2,4,6\}}$. It was further noted in \cite[Example~4.1]{BrandenDLeon} that by \cite[Example~11.5]{HPPlong} the polynomial
    \begin{equation*}
        h_{F_7}+\mu x_2x_4x_6
    \end{equation*}
    is stable only for $\mu=4$. As a side note we would like to mention that this implies in particular that
    \begin{equation*}
        \sum_{S\in\binom{[7]}{3}}|\det(A[S])|\cdot\prod_{i\in S}x_i
    \end{equation*}
    is not stable, giving a negative answer to the question raised in \cite[Remark~4.2]{purbhoo}.
\end{ex}

Now we are ready to disprove the following conjecture by Br{\"a}nd{\'e}n--{D}'Le{\'o}n.

\begin{con}[{\cite[Conjecture~4.2]{BrandenDLeon}}]\label{con:relax}
    Suppose that $M$ has the weak half-plane property. Then so does any relaxation of $M$.
\end{con}

Our counterexample is a suitable relaxation of $P_8$.

\begin{ex}\label{ex:counterrelax}
    Consider again the matroid $P_8$ from \Cref{ex:p8}. The set $X:=\{3,5,6,7\}$ is a circuit-hyperplane of $P_8$. Following \cite{ninematroids} we denote by $P_1$ the relaxation of $P_8$ by $X$. We will prove that $P_1$ does not have the weak half-plane property, although it is a relaxation of $P_8$.    
    Using the {\texttt{Macaulay2}} \cite{M2} package ``Matroids'' \cite{MatroidsArticle} we compute that $\dim(V_{P_1})=10$. Thus by \Cref{ex:p8} and \Cref{rem:embed} the vectors $\delta_X$ and $\iota(u(A))$ span a linear complement of $W_{P_1}$ in $V_{P_1}$. We define $v=\frac{1}{\log(2)}\iota(u(A))$ --- this is just for convenience to get a vector with entries in $\{0,1,2\}$. If $P_1$ has the weak half-plane property, then by \Cref{lem:rescale} there are $a,b>0$ such that
    \begin{equation*}
        F_{a,b}=\sum_{B\text{ basis of }P_8}b^{v_B}x^B+ax_3x_5x_6x_7
    \end{equation*}
    is stable. Now consider the $(0,1)$-th Rayleigh difference
    \begin{equation*}
      \Delta_{0,1}F_{a,b}:=\partial_{x_0}F_{a,b}\partial_{x_1}F_{a,b}-F_{a,b}\partial_{x_0}\partial_{x_1}F_{a,b}.  
    \end{equation*}
    By \cite[Theorem~5.10]{BrandenHPP} the Rayleigh difference $\Delta_{0,1}F_{a,b}$ is globally nonnegative if $F_{a,b}$ is stable. Plugging in $(1,1,t,-1,-1,t)$ for $(x_2,\ldots,x_8)$ yields
    \begin{equation*}
      (\Delta_{0,1}F_{a,b})(1,1,t,-1,-1,t)=-abt^3+(-ab-4b^2+2a+12b+16)t^2+at.  
    \end{equation*}
     Because $a,b>0$, this takes negative values for large enough $t$. Therefore, for no choice of $a,b>0$ the polynomial $F_{a,b}$ is stable.
    
\end{ex}

\begin{rem}
    According to \cite[Proposition~4]{ninematroids} the matroid $P_1$ is not representable over any field. It does not have the half-plane property because it has the matroid $F_7^{-3}$ as a minor which does not have the half-plane property by \cite[Example~11.7]{HPPlong}. Moreover, it does not seem to be representable even in the more general context studied in \cite{skewpartial}. This made $P_1$ a good candidate for being a counterexample to \Cref{con:relax} as representations over certain algebras (that are not necessarily fields) sometimes can still be used to prove the weak half-plane property. See for instance \cite{NonPappus} where the weak half-plane property was proved for the non-Pappus and  the non-Desargues matroid.
\end{rem}

\begin{figure}
\begin{tikzpicture}[line cap=round,line join=round,>=triangle 45,x=0.75cm,y=0.75cm]
\clip(-4,-4) rectangle (4,4);
\draw [line width=2pt] (0,0)-- (-2.6,-1.5);
\draw [line width=2pt] (0,0)-- (2.6,-1.5);
\draw [line width=2pt] (0,3)-- (-2.6,-1.5);
\draw [line width=2pt] (0,3)-- (2.6,-1.5);
\draw [shift={(-4.49,1.5)},line width=2pt]  plot[domain=-0.32:0.32,variable=\t]({1*4.73*cos(\t r)+0*4.73*sin(\t r)},{0*4.73*cos(\t r)+1*4.73*sin(\t r)});
\draw [shift={(4.5,1.5)},line width=2pt]  plot[domain=2.82:3.46,variable=\t]({1*4.74*cos(\t r)+0*4.74*sin(\t r)},{0*4.74*cos(\t r)+1*4.74*sin(\t r)});
\draw [line width=2pt] (-2.6,-1.5)-- (0,-1.5);
\draw [line width=2pt] (0,-1.5)-- (2.6,-1.5);
\fill [color=black] (0,0) circle (4.5pt);
\fill [color=black] (0,3) circle (4.5pt);
\fill [color=black] (-2.6,-1.5) circle (4.5pt);
\fill [color=black] (2.6,-1.5) circle (4.5pt);
\draw[color=black] (-1.3,-0.5) node {\textbf{6}};
\draw[color=black] (1.3,-0.5) node {\textbf{1}};
\draw[color=black] (-1.5,0.9) node {\textbf{5}};
\draw[color=black] (1.5,0.9) node {\textbf{2}};
\fill [color=black] (0,-1.5) circle (4.5pt);
\draw[color=black] (0.5,1.4) node {\textbf{4}};
\draw[color=black] (-0.5,1.4) node {\textbf{7}};
\draw[color=black] (-1.4,-1.7) node {\textbf{0}};
\draw[color=black] (1.4,-1.7) node {\textbf{3}};
\end{tikzpicture}
\begin{tikzpicture}[line cap=round,line join=round,>=triangle 45,x=0.75cm,y=0.75cm]
\clip(-4,-4) rectangle (4,4);
\draw [line width=2pt] (0,3)-- (-3,0);
\draw [line width=2pt] (-3,0)-- (0,-3);
\draw [line width=2pt] (0,-3)-- (3,0);
\draw [line width=2pt] (3,0)-- (0,3);
\draw [shift={(-1.5,4.5)},line width=2pt]  plot[domain=4.39:5.03,variable=\t]({1*4.75*cos(\t r)+0*4.75*sin(\t r)},{0*4.75*cos(\t r)+1*4.75*sin(\t r)});
\draw [shift={(-1.5,-4.5)},line width=2pt]  plot[domain=1.25:1.89,variable=\t]({1*4.74*cos(\t r)+0*4.74*sin(\t r)},{0*4.74*cos(\t r)+1*4.74*sin(\t r)});
\draw [shift={(1.5,-4.5)},line width=2pt]  plot[domain=1.25:1.89,variable=\t]({1*4.74*cos(\t r)+0*4.74*sin(\t r)},{0*4.74*cos(\t r)+1*4.74*sin(\t r)});
\draw [shift={(1.5,4.5)},line width=2pt]  plot[domain=4.39:5.03,variable=\t]({1*4.74*cos(\t r)+0*4.74*sin(\t r)},{0*4.74*cos(\t r)+1*4.74*sin(\t r)});
\fill [color=black] (0,3) circle (4.5pt);
\fill [color=black] (-3,0) circle (4.5pt);
\fill [color=black] (0,-3) circle (4.5pt);
\fill [color=black] (3,0) circle (4.5pt);
\fill [color=black] (0,0) circle (4.5pt);
\end{tikzpicture}
        \caption{The graphs $G_1$ (left) and $G_2$ (right) from \Cref{rem:subd}.}
        \label{fig:graphs}
 \end{figure}
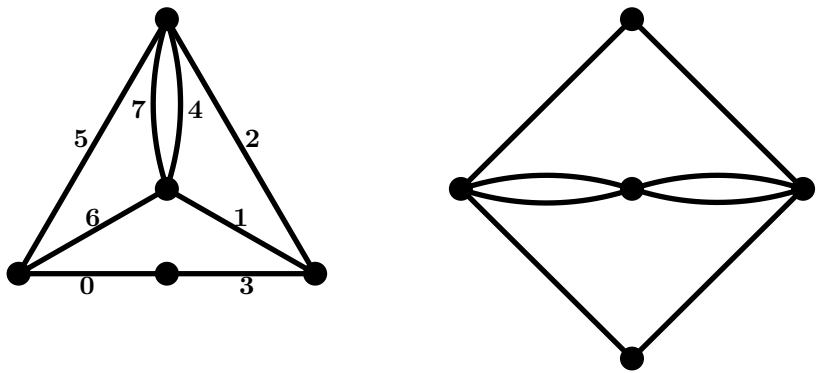

\begin{rem}\label{rem:subd}
    We use the notation as in \Cref{ex:counterrelax}.
    The polynomial $F_{0,1}:=\lim_{a\to0}F_{a,1}$ is the basis generating polynomial of $P_8$. The polynomial $F_{0,0}:=\lim_{a,b\to0}F_{a,b}$ is the basis generating polynomial of the graphical matroid $M(G_1)$ where $G_1$ is depicted in \Cref{fig:graphs} on the left. In total, the regular subdivision of the matroid polytope of $P_8$ defined by $v_B$ has six maximal cells, all of which are matroid polytopes themselves. Four of the corresponding matroids are isomorphic to $M(G_1)$, the two remaining ones are isomorphic to $M(G_2)$ where $G_2$ is the graph on the right of \Cref{fig:graphs}.
    Unfortunately, this knowledge did not help us representing $F_{a,b}$ in a simple way.
\end{rem}

\section{Quaternionic unimodular matroids}\label{sec:qu}
We recall the definition of quaternionic unimodular (QU) matroids. To this end, let $\HH$ denote the skew field of quaternions.

\begin{Def}
    Let $E$ be a finite set and $C\subseteq\HH^E$ a submodule of the free left $\HH$-module $\HH^E$. A nonzero element $x\in C$ is called an \emph{elementary chain} of $C$ if $C$ does not contain a nonzero element whose support is strictly contained in the support of $x$. The submodule $C$ is called \emph{unimodular} if for every elementary chain $x$ of $C$ all nonzero entries of $x$ have the same norm.
\end{Def}

\begin{thm}[{\cite[Theorem~3.7]{skewpartial}}]
    Let $E$ be a finite set and $C\subseteq\HH^E$ unimodular. The set of supports of elementary chains in $C$ is the set of cocircuits of a matroid $M(C)$ on $E$.
\end{thm}

\begin{Def}\label{def:qu}
    A matroid on a finite set $E$ of the form $M(C)$ for some unimodular $C\subseteq\HH^E$ is called \emph{quaternionic unimodular (QU)}.
\end{Def}

The goal of this section is to prove that every QU matroid has the half-plane property. This has been conjectured in \cite[Conjecture~6.9]{skewpartial}.

\begin{Def}[{\cite[page 219]{skewpartial}}]
    Denote by $\varphi\colon\HH\to\C^{2\times 2}$ the map 
\begin{equation*}
    \varphi(a+bi+cj+dk)=
    \begin{pmatrix}
        a+bi& c+di\\
        -c+di& a-bi
    \end{pmatrix}.
\end{equation*}
    We extend $\varphi$ to matrices by applying $\varphi$ entry-wise. Thus for a matrix $A\in\HH^{n\times n}$ we obtain the $2n\times2n$ complex matrix $\varphi(A)$. We define
    \begin{equation*}
      \delta(A):=\sqrt{|\det(\varphi(A))|}.  
    \end{equation*}
\end{Def}

\begin{rem}\label{rem:deltaHom}We collect some basic properties of $\delta$.
\begin{enumerate}
    \item For columns $a_1,\ldots,a_n\in\HH^n$ and $\lambda\in\R$, we have
    \begin{equation*}
      \delta(\lambda a_1\ a_2\ldots a_n)=|\lambda|\cdot\delta(a_1\ldots a_n).  
    \end{equation*}
    This follows directly from the definition.
    \item For matrices $A,B\in\HH^{n\times n}$ we have $\delta(AB)=\delta(A)\delta(B)$ and $\delta(A)=\delta(A^*)$. This is \cite[Lemma~5.2]{skewpartial}.
\end{enumerate}
\end{rem}

For a matrix $A\in\HH^{m\times n}$ with $m\leq n$ and $B\subseteq[n]$ of size $m$ we denote by $A[B]$ the $m\times m$ submatrix of $A$ consisting of the columns indexed by $B$.
The following is a version of the Cauchy--Binet theorem over $\HH$.

\begin{prop}[{\cite[Theorem 5.1]{skewpartial}}]\label{prop:CauchyBinetQuat}
    Let $A\in\HH^{m\times n}$ be a matrix over the quaternions and $m\leq n$. Then
    \begin{equation*}
        \delta(AA^*)=\sum_{\substack{B\subseteq[n]\\|B|=m}}\delta(A[B]A[B]^*).
    \end{equation*}
\end{prop}

\begin{lem}\label{lem:basisIffDelta1}
    Let $E$ be a finite set and $C\subseteq\HH^E$ unimodular. Let $d$ be the rank of the matroid $M=M(C)$.
    There is a $d\times |E|$ matrix $A$ over $\HH$ whose rows form a basis of $C$ such that for all $B\subseteq E$ of size $d$ we have $\delta(A[B])=1$ if $B$ is a basis of $M$ and $\delta(A[B])=0$ otherwise.
\end{lem}

\begin{proof}
 Let $A$ be a matrix over $\HH$ whose rows form a basis of $C$. By \cite[Lemma~3.14]{skewpartial} it has $d$ rows, so $A\in\HH^{d\times |E|}$. After multiplying $A$ from the left by an invertible $d\times d$ matrix over $\HH$, we can assume by \cite[Corollary~3.26]{skewpartial} that $A$ is a \emph{strong QU matrix} in the sense of \cite[Definition~3.23]{skewpartial}. Now the statement of the lemma follows from {\cite[Claim 5.4.1]{skewpartial}}.
\end{proof}

\begin{thm}[{\cite[Conjecture~6.9]{skewpartial}}]
    Every QU matroid has the half-plane property.
\end{thm}
\begin{proof}
  Let $E=\{1,\ldots,m\}$ and $C\subseteq\HH^m$ be unimodular. Denote by $M=M(C)$ the associated QU matroid and let $A\in\HH^{d\times m}$ be a matrix as in \Cref{lem:basisIffDelta1}. We have to show that the bases generating polynomial $h_M$ of $M$ is stable. This is equivalent to $h_M^2$ being stable. We prove this by showing that $h_M^2$ agrees with the stable polynomial
  \begin{equation*}
    \det\left(\varphi(A)\begin{pmatrix}x_1&&&&\\&x_1&&&\\&&\ddots&&\\&&&x_m&\\&&&&x_m\end{pmatrix}\varphi(A)^*\right).  
  \end{equation*}
    It suffices to show that these two polynomials agree on the positive orthant. We denote by $a_1,\ldots,a_m\in\HH^d$ the columns of $A$. Let $x\in\R_{>0}^m$ and write $x=(x_1^2,\ldots,x_m^2)$ for $x_1,\ldots,x_m\in\R_{>0}$. Then we have
    \allowdisplaybreaks\begin{align}
        h_M(x)^2&=\left(\sum_{\substack{B\subseteq E\\|B|=d}}\delta(A[B]A[B]^*)x^B\right)^2\\
        &=\left(\sum_{\substack{B\subseteq E\\|B|=d}}\delta((x_1a_1\ldots x_ma_m)[B](x_1a_1\ldots x_ma_m)[B]^*)\right)^2\\
        &=\delta((x_1a_1\ldots x_ma_m)(x_1a_1\ldots x_ma_m)^*)^2\\
        &=|\det(\varphi(x_1a_1\ldots x_ma_m)\varphi(x_1a_1\ldots x_ma_m)^*)|\\
        &=|\det\left(\varphi(A)\begin{pmatrix}x_1^2&&&&\\&x_1^2&&&\\
        &&\ddots&&\\&&&x_m^2&\\&&&&x_m^2\end{pmatrix}\varphi(A)^*\right)|\\
        &=\det\left(\varphi(A)\begin{pmatrix}x_1^2&&&&\\&x_1^2&&&\\
        &&\ddots&&\\&&&x_m^2&\\&&&&x_m^2\end{pmatrix}\varphi(A)^*\right).
    \end{align}
    Here we have equality in (2) by \Cref{rem:deltaHom} and \Cref{lem:basisIffDelta1}. (3) holds  by \Cref{rem:deltaHom}. For (4) we use \Cref{prop:CauchyBinetQuat} and (5) follows from the definition of $\delta$. (6) and (7) are obvious.
\end{proof}

\bigskip

%  \noindent \textbf{Acknowledgements.}
% We would like to thank somebody.

\bibliographystyle{alpha}
\bibliography{biblio}
 \end{document}

%% file: operators.tex
\usepackage{mathrsfs}
\newcommand{\Ima}{\operatorname{Im}}

\newcommand{\supp}{\operatorname{supp}}

\newcommand{\rank}{\operatorname{rank}} %rank of a free module

%% file: fonts.tex
\newcommand{\cB}{{\mathcal B}}

\newcommand{\C}{{\mathbb C}}
\newcommand{\HH}{{\mathbb H}}
\newcommand{\R}{{\mathbb R}}

\newcommand{\N}{{\mathbb N}}